\documentclass[11pt]{article}
\usepackage{geometry}
\usepackage{amsmath,amsthm,amssymb,enumerate,latexsym}
\usepackage[colorlinks=true,linkcolor=blue,citecolor=blue,urlcolor=black,pdfpagelabels=false]{hyperref}

\usepackage{array}

\usepackage{fix-cm}

\newtheorem{theorem}{Theorem}
\newtheorem{lemma}[theorem]{Lemma}
\newtheorem{conjecture}[theorem]{Conjecture}

\theoremstyle{definition}
\newtheorem{example}[theorem]{Example}
\newtheorem*{exercise*}{Exercise}

\newenvironment{acknowledgment}{\bigskip\textbf{Acknowledgements.}}{}
\newcommand{\email}[1]{{\textit{Email:} \texttt{#1}}}

\makeatletter
\newcommand{\dedicatory}[1]{%
  \let\@@@oldtitle\@title%
  \gdef\@title{\@@@oldtitle{\\[0.5em]{\hfill\small\sl #1}\\[0.75em]}}%
}
\makeatother

\newcommand{\assign}{:=}
\newcommand{\md}{\mathrm{d}}

\begin{document}

\title{Ap\'ery Limits: Experiments and Proofs}

\dedicatory{Dedicated to the memory of Jon Borwein}

\author{
  Marc Chamberland
  \thanks{Grinnell College, \email{chamberl@math.grinnell.edu}}
  \and
  Armin Straub
  \thanks{University of South Alabama, \email{straub@southalabama.edu}}
}

\date{November 6, 2020}

\maketitle

\begin{abstract}
  An important component of Ap\'ery's proof that $\zeta (3)$ is irrational
  involves representing $\zeta (3)$ as the limit of the quotient of two
  rational solutions to a three-term recurrence. We present various approaches
  to such Ap\'ery limits and highlight connections to continued fractions as
  well as the famous theorems of Poincar\'e and Perron on difference
  equations. In the spirit of Jon Borwein, we advertise an
  experimental-mathematics approach by first exploring in detail a simple but
  instructive motivating example. We conclude with various open problems.
\end{abstract}

\section{Introduction}

A fundamental ingredient of Ap\'ery's groundbreaking proof \cite{apery} of
the irrationality of $\zeta (3)$ is the binomial sum
\begin{equation}
  A (n) = \sum_{k = 0}^n \binom{n}{k}^2 \binom{n + k}{k}^2 \label{eq:apery3}
\end{equation}
and the fact that it satisfies the three-term recurrence
\begin{equation}
  (n + 1)^3 u_{n + 1} = (2 n + 1) (17 n^2 + 17 n + 5) u_n - n^3 u_{n - 1}
  \label{eq:apery3:rec}
\end{equation}
with initial conditions $A (0) = 1$, $A (1) = 5$ --- or, equivalently but more
naturally, $A (- 1) = 0$, $A (0) = 1$. Now let $B (n)$ be the solution to
{\eqref{eq:apery3:rec}} with $B (0) = 0$ and $B (1) = 1$. Ap\'ery showed
that
\begin{equation}
  \lim_{n \rightarrow \infty} \frac{B (n)}{A (n)} = \frac{\zeta (3)}{6}
  \label{eq:apery3:lim}
\end{equation}
and that the rational approximations resulting from the left-hand side
converge too rapidly to $\zeta (3)$ for $\zeta (3)$ itself to be rational. For
details, we recommend the engaging account \cite{alf} of Ap\'ery's proof.
In the sequel, we will not pursue questions of irrationality further. Instead,
our focus will be on limits, like {\eqref{eq:apery3:lim}}, of quotients of
solutions to linear recurrences. Such limits are often called {\emph{Ap\'ery
limits}} \cite{avz-apery-limits}, \cite{yang-apery-limits}.

Jon Borwein was a tireless advocate and champion of experimental mathematics
and applied it with fantastic success. Jon was also a pioneer of teaching
experimental mathematics, whether through numerous books, such as
\cite{bb-exp1}, or in the classroom (the second author is grateful for the
opportunity to benefit from both). Before collecting known results on
Ap\'ery limits and general principles, we therefore find it fitting to
explore in detail, in Section~\ref{sec:delannoy}, a simple but instructive
example using an experimental approach. We demonstrate how to discover the
desired Ap\'ery limit; and we show, even more importantly, how the
exploratory process naturally leads us to discover additional structure that
is helpful in understanding this and other such limits. We hope that the
detailed discussion in Section~\ref{sec:delannoy} may be particularly useful
to those seeking to integrate experimental mathematics into their own
teaching.

After suggesting further examples in Section~\ref{sec:search}, we explain the
observations made in Section~\ref{sec:delannoy} by giving in
Section~\ref{sec:diffeq} some background on difference equations, introducing
the Casoratian and the theorems of Poincar\'e and Perron. In
Section~\ref{sec:cf}, we connect with continued fractions and observe that,
accordingly translated, many of the simpler examples are instances of
classical results in the theory of continued fractions. We then outline in
Section~\ref{sec:pf} several methods used in the literature to establish
Ap\'ery limits. To illustrate the limitations of these approaches, we
conclude with several open problems in Sections~\ref{sec:franel:d} and
\ref{sec:open}.

Creative telescoping --- including, for instance, Zeilberger's algorithm and
the Wilf--Zeilberger (WZ) method --- refers to a powerful set of tools that,
among other applications, allow us to algorithmically derive the recurrence
equations, like {\eqref{eq:apery3:rec}}, that are satisfied by a given sum,
like {\eqref{eq:apery3}}. In fact, as described in \cite{alf}, Zagier's
proof of Ap\'ery's claim that the sums {\eqref{eq:apery3}} and
{\eqref{eq:apery3:2}} both satisfy the recurrence {\eqref{eq:apery3:rec}} may
be viewed as giving birth to the modern method of creative telescoping. For an
excellent introduction we refer to \cite{aeqb}. In the sequel, all claims
that certain sums satisfy a recurrence can be established using creative
telescoping.

\section{A motivating example}\label{sec:delannoy}

At the end of van der Poorten's account \cite{alf} of Ap\'ery's proof, the
reader is tasked with the exercise to consider the sequence
\begin{equation}
  A (n) = \sum_{k = 0}^n \binom{n}{k} \binom{n + k}{k} \label{eq:delannoy}
\end{equation}
and to apply to it Ap\'ery's ideas to conclude the irrationality of $\ln
(2)$. In this section, we will explore this exercise with an experimental
mindset but without using the general tools and connections described later in
the paper. In particular, we hope that the outline below could be handed out
in an undergraduate experimental-math class and that the students could (with
some help, depending on their familiarity with computer algebra systems)
reproduce the steps, feel intrigued by the observations along the way, and
then apply by themselves a similar approach to explore variations or
extensions of this exercise. Readers familiar with the topic may want to skip
ahead.

The numbers {\eqref{eq:delannoy}} are known as the central Delannoy numbers
and count lattice paths from $(0, 0)$ to $(n, n)$ using the steps $(0, 1)$,
$(1, 0)$ and $(1, 1)$. They satisfy the recurrence
\begin{equation}
  (n + 1) u_{n + 1} = 3 (2 n + 1) u_n - n u_{n - 1} \label{eq:delannoy:rec}
\end{equation}
with initial conditions $A (- 1) = 0$, $A (0) = 1$. Now let $B (n)$ be the
sequence satisfying the same recurrence with initial conditions $B (0) = 0$,
$B (1) = 1$. Numerically, we observe that the quotients $Q (n) = B (n) / A
(n)$,
\begin{equation*}
  (Q (n))_{n \geq 0} = \left(0, \frac{1}{3}, \frac{9}{26},
   \frac{131}{378}, \frac{445}{1284}, \frac{34997}{100980},
   \frac{62307}{179780}, \frac{2359979}{6809460}, \ldots \right),
\end{equation*}
appear to converge rather quickly to a limit
\begin{equation*}
  L \assign \lim_{n \rightarrow \infty} Q (n) = 0.34657359 \ldots
\end{equation*}
When we try to estimate the speed of convergence by computing the difference
$Q (n) - Q (n - 1)$ of consecutive terms, we find
\begin{equation*}
  (Q (n) - Q (n - 1))_{n \geq 1} = \left(\frac{1}{3}, \frac{1}{78},
   \frac{1}{2457}, \frac{1}{80892}, \frac{1}{2701215}, \frac{1}{90770922},
   \ldots \right) .
\end{equation*}
This suggests the probably-unexpected fact that these are all reciprocals of
integers. Something interesting must be going on here! However, a cursory
look-up of the denominators in the {\emph{\textit{On-Line Encyclopedia of
Integer Sequences}}} (OEIS) \cite{oeis} does not result in a match. (Were we
to investigate the factorizations of these integers, we might at this point
discover the case $x = 1$ of {\eqref{eq:delannoy:x:Qdiff}}. But we hold off on
exploring that observation and continue to focus on the speed of convergence.)
By, say, plotting the logarithm of $Q (n) - Q (n - 1)$ versus $n$, we are led
to realize that the number of digits to which $Q (n - 1)$ and $Q (n)$ agree
appears to increase (almost perfectly) linearly. This means that $Q (n)$
converges to $L$ exponentially.

\begin{exercise*}
  For a computational challenge, quantify the exponential convergence by
  conjecturing an exact value for the limit of $(Q (n + 1) - Q (n)) / (Q (n) -
  Q (n - 1))$ as $n \rightarrow \infty$. Then connect that value to the
  recurrence {\eqref{eq:delannoy:rec}}.
\end{exercise*}

At this point, we are confident that, say,
\begin{equation}
  Q (50) = 0.34657359027997265470861606072908828403775006718 \ldots
  \label{eq:delannoy:Q50}
\end{equation}
agrees with the limit $L$ to more than $75$ digits. The ability to recognize
constants from numerical data is a powerful asset in an experimental
mathematician's toolbox. Several approaches to constant recognition are
lucidly described in \cite[Section~6.3]{bb-exp1}. The crucial ingredients
are integer relation algorithms such as PSLQ or those based on lattice
reduction algorithms like LLL. Readers new to constant recognition may find
the {\emph{Inverse Symbolic Calculator}} of particular value. This web
service, created by Jon Borwein, Peter Borwein and Simon Plouffe, automates
the constant-recognition process: it asks for a numerical approximation as
input and determines, if successful, a suggested exact value. For instance,
given {\eqref{eq:delannoy:Q50}}, it suggests that
\begin{equation*}
  L = \frac{1}{2} \ln (2),
\end{equation*}
which one can then easily confirm further to any desired precision. Of course,
while this provides overwhelming evidence, it does not constitute a proof.
Given the success of our exploration, a natural next step would be to repeat
this inquiry for the sequence of polynomials
\begin{equation}
  A_x (n) = \sum_{k = 0}^n \binom{n}{k} \binom{n + k}{k} x^k,
  \label{eq:delannoy:x}
\end{equation}
which satisfies the recurrence {\eqref{eq:delannoy:rec}} with the term $3 (2 n
+ 1)$ replaced by $(2 x + 1) (2 n + 1)$. An important principle to keep in
mind here is that introducing an additional parameter, like the $x$ in
{\eqref{eq:delannoy:x}}, can make the underlying structure more apparent; and
this may be crucial both for guessing patterns and for proving our assertions.
Now define the secondary solution $B_x (n)$ satisfying the recurrence with
$B_x (0) = 0$, $B_x (1) = 1$. Then, if we compute the difference of quotients
$Q_x (n) = B_x (n) / A_x (n)$ as before, we find that
\begin{equation*}
  (Q_x (n) - Q_x (n - 1))_{n \geq 1} = \left(\frac{1}{1 + 2 x},
   \frac{1}{2 (1 + 2 x) (1 + 6 x + 6 x^2)}, \ldots \right) .
\end{equation*}
Extending our earlier observation, these now appear to be the reciprocals of
polynomials with integer coefficients. Moreover, in factored form, we are
immediately led to conjecture that
\begin{equation}
  Q_x (n) - Q_x (n - 1) = \frac{1}{n A_x (n) A_x (n - 1)} .
  \label{eq:delannoy:x:Qdiff}
\end{equation}
Note that, since $Q_x (0) = 0$, this implies
\begin{equation}
  Q_x (N) = \sum_{n = 1}^N (Q_x (n) - Q_x (n - 1)) = \sum_{n = 1}^N \frac{1}{n
  A_x (n) A_x (n - 1)} \label{eq:delannoy:x:Qsum}
\end{equation}
and hence provides another way to compute the limit $L_x = \lim_{n \rightarrow
\infty} Q_x (n)$ as
\begin{equation*}
  L_x = \sum_{n = 1}^{\infty} \frac{1}{n A_x (n) A_x (n - 1)},
\end{equation*}
which avoids reference to the secondary solution $B_x (n)$.

Can we experimentally identify the limit $L_x$? One approach could be to
select special values for $x$ and then proceed as we did for $x = 1$. For
instance, we might numerically compute and then identify the following values:
\begin{equation*}
  \renewcommand{\arraystretch}{1.3}
  \begin{array}{|l|l|l|l|}
     \hline
     & x = 1 & x = 2 & x = 3\\
     \hline
     L_x & \tfrac{1}{2} \ln (2) & \tfrac{1}{2} \ln \left(\tfrac{3}{2} \right)
     & \tfrac{1}{2} \ln \left(\tfrac{4}{3} \right)\\
     \hline
   \end{array}
\end{equation*}
We are lucky and the emerging pattern is transparent, suggesting that
\begin{equation}
  L_x = \frac{1}{2} \ln \left(1 + \frac{1}{x} \right) .
  \label{eq:delannoy:x:L}
\end{equation}
A possibly more robust approach to identifying $L_x$ empirically is to fix
some values of $n$ and then expand the $Q_x (n)$, which are rational functions
in $x$, into power series. If the initial terms of these power series appear
to converge as $n \rightarrow \infty$ to identifiable values, then it is
reasonable to expect that these values are the initial terms of the power
series for the limit $L_x$. However, expanding around $x = 0$, we quickly
realize that the power series
\begin{equation*}
  Q_x (n) = \sum_{k = 0}^{\infty} q_k^{(n)} x^k
\end{equation*}
do not stabilize as $n \rightarrow \infty$, but that the coefficients increase
in size: for instance, we find empirically that
\begin{equation*}
  q_1^{(n)} = - n (n + 1), \quad q_2^{(n)} = \frac{1}{8} n (n + 1) (5 n^2 + 5
   n + 6),
\end{equation*}
and it appears that, for $k \geq 1$, $q_k^{(n)}$ is a polynomial in $n$
of degree $2 k$. Expanding the $Q_x (n)$ instead around some nonzero value of
$x$ --- say, $x = 1$ --- is more promising. Writing
\begin{equation*}
  Q_x (n) = \sum_{k = 0}^{\infty} r_k^{(n)} (x - 1)^k,
\end{equation*}
we observe empirically that
\begin{equation*}
  \left( \lim_{n \rightarrow \infty} r_k^{(n)} \right)_{k \geq 1} = \left(-
   \frac{1}{4}, \frac{3}{16}, - \frac{7}{48}, \frac{15}{128}, \ldots \right) .
\end{equation*}
Once we realize that the denominators are multiples of $k$, it is not
difficult to conjecture that
\begin{equation}
  \lim_{n \rightarrow \infty} r_k^{(n)} = (- 1)^k \frac{2^k - 1}{k \cdot 2^{k
  + 1}} \label{eq:delannoy:x:Q:c1}
\end{equation}
for $k \geq 1$. From our initial exploration, we already know that
$\lim_{n \rightarrow \infty} r_0^{(n)} = \frac{1}{2} \ln (2)$ but we could
also have (re)guessed this value as the formal limit of the right-hand side of
{\eqref{eq:delannoy:x:Q:c1}} as $k \rightarrow 0$ (forgetting that $k$ is
really an integer). Anyway, {\eqref{eq:delannoy:x:Q:c1}} suggests that
\begin{equation*}
  L_x = L_1 + \sum_{k = 1}^{\infty} (- 1)^k \frac{2^k - 1}{k \cdot 2^{k + 1}}
   \, (x - 1)^k = \frac{1}{2} \ln (2) + \frac{1}{2} \ln \left(\frac{x + 1}{2
   x} \right),
\end{equation*}
leading again to {\eqref{eq:delannoy:x:L}}. Finally, our life is easiest if we
look at the power series of $Q_x (n)$ expanded around $x = \infty$. In that
case, we find that the power series of $Q_x (n)$ and $Q_x (n + 1)$ actually
agree through order $x^{- 2 n}$. In hindsight --- and to expand our
experimental reach, it is always a good idea to reflect on the new data in
front of us --- this is a consequence of {\eqref{eq:delannoy:x:Qdiff}} and the
fact that $A_x (n)$ has degree $n$ in $x$ (so that $A_x (n) A_x (n - 1)$ has
degree $2 n - 1$). Therefore, from just the case $n = 3$ we are confident that
\begin{equation*}
  L_x = Q_x (3) + O (x^{- 7}) = \frac{1}{2 x} - \frac{1}{4 x^2} + \frac{1}{6
   x^3} - \frac{1}{8 x^4} + \frac{1}{10 x^5} - \frac{1}{12 x^6} + O (x^{- 7})
   .
\end{equation*}
At this point the pattern is evident, and we arrive, once more, at the
conjectured formula {\eqref{eq:delannoy:x:L}} for $L_x$.

\section{Searching for Ap\'ery limits}\label{sec:search}

Inspired by the approach laid out in the previous section, one can search for
other Ap\'ery limits as follows:
\begin{enumerate}
  \item Pick a binomial sum $A (n)$ and, using creative telescoping, compute a
  recurrence satisfied by $A (n)$.
  
  \item Compute the initial terms of a secondary solution $B (n)$ to the
  recurrence.
  
  \item Try to identify $\lim_{n \rightarrow \infty} B (n) / A (n)$ (either
  numerically or as a power series in an additional parameter).
\end{enumerate}
It is important to realize, as will be indicated in Section~\ref{sec:cf}, that
if the binomial sum $A (n)$ satisfies a three-term recurrence, then the
Ap\'ery limit can be expressed as a continued fraction and compared to the
(rather staggering) body of known results \cite{wall-contfrac},
\cite{jt-contfrac}, \cite{handbook-contfrac}, \cite{bvsz-cf}.

Of course, the final step is to prove and/or generalize those discovered
results that are sufficiently appealing. One benefit of an experimental
approach is that we can discover results, connections and generalizations, as
well as discard less-fruitful avenues, before (or while!) working out a
complete proof. Ideally, the processes of discovery and proof inform each
other at every stage. For instance, experimentally finding a generalization
may well lead to a simplified proof, while understanding a small piece of a
puzzle can help set the direction of follow-up experiments. Jon Borwein's
extensive legacy is filled with such delightful examples.

Of course, one could just start with a recurrence; however, selecting a
binomial sum increases the odds that the recurrence has desirable properties
(it is a difficult open problem to ``invert creative telescoping'' in the
sense of producing a binomial sum satisfying a given recurrence). Some simple
suggestions for binomial sums, as well as the corresponding Ap\'ery limits,
are as follows (in each case, we choose the secondary solution with initial
conditions $B (0) = 0$, $B (1) = 1$).
\begin{equation*}
  \setlength{\extrarowheight}{7pt}
  \renewcommand{\arraystretch}{1.3}
  \begin{array}{|>{\displaystyle}l|>{\displaystyle}l|}
     \hline
     \sum_{k = 0}^n \binom{n}{2 k} x^k & \frac{1}{\sqrt{x}} \quad
     \text{(around $x = 1$)}\\
     \hline
     \sum_{k = 0}^n \binom{n - k}{k} x^k & \frac{2}{1 + \sqrt{1 + 4 x}} \quad
     \text{(around $x = 0$)}\\
     \hline
     \sum_{k = 0}^n \binom{n}{k} \binom{n - k}{k} x^k & \frac{\arctan \left(\sqrt{4 x - 1} \right)}{\sqrt{4 x - 1}} \quad \text{(around $x =
     \tfrac{1}{4}$)}\\
     \hline
   \end{array}
\end{equation*}
\begin{example}
  \label{eg:arctan}Setting $x = 1 / 2$ in the last instance leads to the limit
  being $\arctan (1) = \pi / 4$ and therefore to a way of computing $\pi$ as
  \begin{equation*}
    \pi = \lim_{n \rightarrow \infty} \frac{4 B (n)}{A (n)},
  \end{equation*}
  where $A (n)$ and $B (n)$ both solve the recurrence $(n + 1) u_{n + 1} = (2
  n + 1) u_n + n u_{n - 1}$ with $A (- 1) = 0$, $A (0) = 1$ and $B (0) = 0$,
  $B (1) = 1$. In an experimental-math class, this could be used to segue into
  the fascinating world of computing $\pi$, a topic to which Jon Borwein,
  sometimes admiringly referred to as Dr.~Pi, has contributed so much --- one
  example being the groundbreaking work in \cite{borwein-piagm} with his
  brother Peter. Let us note that this is not a new result. Indeed, with the
  substitution $z = \sqrt{4 x - 1}$, it follows from the discussion in
  Section~\ref{sec:cf} that the Ap\'ery limit in question is equivalent to
  the well-known continued fraction
  \begin{equation*}
    \arctan (z) = \frac{z}{1 +} \, \frac{1^2 z^2}{3 +} \, \frac{2^2 z^2}{5 +}
     \, \cdots \, \frac{n^2 z^2}{(2 n + 1) +} \, \cdots
  \end{equation*}
  \cite[p.~343, eq.~(90.3)]{wall-contfrac}. The reader finds, for instance,
  in \cite[Theorem~2]{bcp-cf-tails} that this continued fraction, as well as
  corresponding ones for the tails of $\arctan (z)$, is a special case of
  Gauss' famous continued fraction for quotients of hypergeometric functions
  ${}_2 F_1$. We hope that some readers and students, in particular, enjoy the
  fact that they are able to rediscover such results themselves.
\end{example}

\begin{example}
  For more challenging explorations, the reader is invited to consider the
  binomial sums
  \begin{equation*}
    A_x (n) = \sum_{k = 0}^n \binom{n}{k} \binom{n + k}{k}^2 x^k, \quad
     \sum_{k = 0}^n \binom{n}{k} \binom{n + k}{k}^3 x^k,
  \end{equation*}
  and to compare the findings with those by Zudilin
  \cite{zudilin-appr-polylog} who obtains simultaneous approximations to the
  logarithm, dilogarithm and trilogarithm.
\end{example}

\begin{example}
  \label{eg:cy}Increasing the level of difficulty further, one may consider,
  for instance, the binomial sum
  \begin{equation*}
    A (n) = \sum_{k = 0}^n \binom{n}{k}^2 \binom{3 k}{n},
  \end{equation*}
  which is an appealing instance, randomly selected from many others, for
  which Almkvist, van Straten and Zudilin \cite[Section~4,
  \#219]{avz-apery-limits} have numerically identified an Ap\'ery limit (in
  this case, depending on the initial conditions of the secondary solution,
  the Ap\'ery limit can be empirically expressed as a rational multiple of
  $\pi^2$ or of the $L$-function evaluation $L_{- 3} (2)$, or, in general, a
  linear combination of those). To our knowledge, proving the conjectured
  Ap\'ery limits for most cases in \cite[Section~4]{avz-apery-limits},
  including the one above, remains open. While the techniques discussed in
  Section~\ref{sec:pf} can likely be used to prove some individual limits, it
  would be of particular interest to establish all these Ap\'ery limits in a
  uniform fashion.
\end{example}

Choosing an appropriate binomial sum as a starting point, the present approach
could be used to construct challenges for students in an experimental-math
class, with varying levels of difficulty (or that students could explore
themselves with minimal guidance). As illustrated by Example~\ref{eg:arctan},
simple cases can be connected with existing classical results, and
opportunities abound to connect with other topics such as hypergeometric
functions, computer algebra, orthogonal polynomials, or Pad\'e
approximation, which we couldn't properly discuss here. However, much about
Ap\'ery limits is not well understood and we believe that more serious
investigations, possibly along the lines outlined here, can help improve our
understanding. To highlight this point, we present in
Sections~\ref{sec:franel:d} and \ref{sec:open} several specific open problems
and challenges.

\section{Difference equations}\label{sec:diffeq}

In our initial motivating example, we started with a solution $A (n)$ to the
three-term recurrence {\eqref{eq:delannoy:rec}} and considered a second,
linearly independent solution $B (n)$ of that same recurrence. We then
discovered in {\eqref{eq:delannoy:x:Qsum}} that
\begin{equation*}
  B (n) = A (n) \sum_{k = 1}^n \frac{1}{k A (k) A (k - 1)} .
\end{equation*}
That the secondary solution is expressible in terms of the primary solution is
a consequence of a general principle in the theory of difference equations,
which we outline in this section. For a gentle introduction to difference
equations, we refer to \cite{kelley-peterson-diff}.

Consider the general homogeneous linear difference equation
\begin{equation}
  u (n + d) + p_{d - 1} (n) u (n + d - 1) + \cdots + p_1 (n) u (n + 1) + p_0
  (n) u (n) = 0 \label{eq:de}
\end{equation}
of order $d$, where we normalize the leading coefficient to $1$. If $u_1 (n),
\ldots, u_d (n)$ are solutions to {\eqref{eq:de}}, then their
{\emph{Casoratian}} $w (n)$ is defined as
\begin{equation*}
  w (n) = \det \begin{bmatrix}
     u_1 (n) & u_2 (n) & \cdots & u_d (n)\\
     u_1 (n + 1) & u_2 (n + 1) & \cdots & u_d (n + 1)\\
     \vdots & \vdots & \ddots & \vdots\\
     u_1 (n + d - 1) & u_2 (n + d - 1) & \cdots & u_d (n + d - 1)
   \end{bmatrix} .
\end{equation*}
This is the discrete analog of the Wronskian that is discussed in most
introductory courses on differential equations. By applying the difference
equation {\eqref{eq:de}} to the last row in $w (n + 1)$ and then subtracting
off multiples of earlier rows, one finds that the Casoratian satisfies
\cite[Lemma~3.1]{kelley-peterson-diff}
\begin{equation*}
  w (n + 1) = (- 1)^d p_0 (n) w (n)
\end{equation*}
and hence
\begin{equation}
  w (n) = (- 1)^{d n} p_0 (0) p_0 (1) \cdots p_0 (n - 1) w (0) .
  \label{eq:casoratian:rec}
\end{equation}
In the case of second order difference equations ($d = 2$), we have
\begin{equation*}
  \frac{u_2 (n + 1)}{u_1 (n + 1)} - \frac{u_2 (n)}{u_1 (n)} = \frac{u_1 (n)
   u_2 (n + 1) - u_1 (n + 1) u_2 (n)}{u_1 (n) u_1 (n + 1)} = \frac{w (n)}{u_1
   (n) u_1 (n + 1)},
\end{equation*}
which implies that we can construct a second solution from a given solution as
follows.

\begin{lemma}
  Let $d = 2$ and suppose that $u_1 (n)$ solves {\eqref{eq:de}} and that $u_1
  (n) \neq 0$ for all $n \geq 0$. Then a second solution of
  {\eqref{eq:de}} is given by
  \begin{equation}
    u_2 (n) = u_1 (n) \sum_{k = 0}^{n - 1} \frac{w (k)}{u_1 (k) u_1 (k + 1)},
    \label{eq:u2:casoratian}
  \end{equation}
  where $w (k) = p_0 (0) p_0 (1) \cdots p_0 (k - 1)$.
\end{lemma}

Here we have normalized the solution $u_2$ by choosing $w (0) = 1$: this
entails $u_2 (0) = 0$ and $u_2 (1) = 1 / u_1 (0)$. Note also that if $p_0 (0)
\neq 0$, then $w (1) \neq 0$, which implies that the solution $u_2$ is
linearly independent from $u_1$.

\begin{example}
  For the Delannoy difference equation {\eqref{eq:delannoy:rec}} and the
  solutions $A (n)$, $B (n)$ with initial conditions as before, we have $d =
  2$ and $p_0 (n) = (n + 1) / (n + 2)$, hence $w (n) = 1 / (n + 1)$. In
  particular, equation {\eqref{eq:u2:casoratian}} is equivalent to
  {\eqref{eq:delannoy:x:Qsum}}.
\end{example}

Now suppose that $p_k (n) \rightarrow c_k$ as $n \rightarrow \infty$, for each
$k \in \{ 0, 1, \ldots, d - 1 \}$. Then the characteristic polynomial of the
recurrence {\eqref{eq:de}} is, by definition,
\begin{equation*}
  \lambda^d + c_{d - 1} \lambda^{d - 1} + \cdots + c_1 \lambda + c_0 =
   \prod_{k = 1}^d (\lambda - \lambda_k)
\end{equation*}
with characteristic roots $\lambda_1, \ldots, \lambda_d$. Poincare's famous
theorem \cite[Theorem~5.1]{kelley-peterson-diff} states, under a modest
additional hypothesis, that each nontrivial solution to {\eqref{eq:de}} has
asymptotic growth according to one of the characteristic roots.

\begin{theorem}[Poincar\'e]
  \label{thm:poincare}Suppose further that the characteristic roots have
  distinct moduli. If $u (n)$ solves the recurrence {\eqref{eq:de}}, then
  either $u (n) = 0$ for all sufficiently large $n$, or
  \begin{equation}
    \lim_{n \rightarrow \infty} \frac{u (n + 1)}{u (n)} = \lambda_k
    \label{eq:poincare}
  \end{equation}
  for some $k \in \{ 1, \ldots, d \}$.
\end{theorem}

If, in addition, $p_0 (n) \neq 0$ for all $n \geq 0$ (so that, by
{\eqref{eq:casoratian:rec}}, the Casoratian $w (n)$ is either zero for all $n$
or nonzero for all $n$), then Perron's theorem guarantees that, for each $k$,
there exists a solution such that {\eqref{eq:poincare}} holds.

\section{Continued Fractions}\label{sec:cf}

In this section, we briefly connect with (irregular) continued fractions
\begin{equation*}
  C = b_0 + \frac{a_1}{b_1 +} \, \frac{a_2}{b_2 +} \, \frac{a_3}{b_3 +}
   \ldots \assign b_0 + \frac{a_1}{b_1 + \frac{a_2}{b_2 + \frac{a_3}{b_3 +
   \ldots}}},
\end{equation*}
as introduced, for instance, in \cite{jt-contfrac},
\cite[Entry~12.1]{berndtII} or \cite[Chapter~9]{bvsz-cf}. The $n$-th
convergent of $C$ is
\begin{equation*}
  C_n = b_0 + \frac{a_1}{b_1 +} \, \frac{a_2}{b_2 +} \, \ldots \,
   \frac{a_n}{b_n} .
\end{equation*}
It is well known, and readily proved by induction, that $C_n = B (n) / A (n)$
where both $A (n)$ and $B (n)$ solve the second-order recurrence
\begin{equation}
  u_n = b_n u_{n - 1} + a_n u_{n - 2} \label{eq:cf:rec}
\end{equation}
with $A (- 1) = 0$, $A (0) = 1$ and $B (- 1) = 1$, $B (0) = b_0$. (Note that,
if $b_0 = 0$, then the initial conditions for $B (n)$ are equivalent to $B (0)
= 0$, $B (1) = a_1$.)

Conversely, see \cite[Theorem~9.4]{bvsz-cf}, two such solutions to a
second-order difference equation with non-vanishing Casoratian correspond to a
unique continued fraction. In particular, Ap\'ery limits $\lim_{n
\rightarrow \infty} B (n) / A (n)$ arising from second-order difference
equations can be equivalently expressed as continued fractions.

\begin{example}
  The Ap\'ery limit {\eqref{eq:delannoy:x:L}} is equivalent to the continued
  fraction
  \begin{equation*}
    \frac{1}{2} \ln \left(1 + \frac{1}{x} \right) = \frac{1}{(2 x + 1) -}
     \, \frac{1^2}{3 (2 x + 1) -} \, \frac{2^2}{5 (2 x +
     1) -} \cdots \frac{n^2}{(2 n + 1) (2 x + 1) -} \cdots
  \end{equation*}
  \cite[p.~343, eq.~(90.4)]{wall-contfrac}. To see this, it suffices to note
  that, if $A_x (n)$ and $B_x (n)$ are as in Section~\ref{sec:delannoy}, then
  $n!A_x (n)$ and $n!B_x (n)$ solve the adjusted recurrence
  \begin{equation*}
    u_{n + 1} = (2 x + 1) (2 n + 1) u_n - n^2 u_{n - 1}
  \end{equation*}
  of the form {\eqref{eq:cf:rec}}.
\end{example}

The interested reader can find a detailed discussion of the continued
fractions corresponding to Ap\'ery's limits for $\zeta (2)$ and $\zeta (3)$
in \cite[Section~9.5]{bvsz-cf}, which then proceeds to proving the
respective irrationality results.

\section{Proving Ap\'ery limits}\label{sec:pf}

In the sequel, we briefly indicate several approaches towards proving
Ap\'ery limits. In case of the Ap\'ery numbers {\eqref{eq:apery3}},
Ap\'ery established the limit {\eqref{eq:apery3:lim}} by finding the
explicit expression
\begin{equation}
  B (n) = \frac{1}{6} \sum_{k = 0}^n \binom{n}{k}^2 \binom{n + k}{k}^2 \left(\sum_{j = 1}^n \frac{1}{j^3} + \sum_{m = 1}^k \frac{(- 1)^{m - 1}}{2 m^3
  \binom{n}{m} \binom{n + m}{m}} \right) \label{eq:apery3:2}
\end{equation}
for the secondary solution $B (n)$. Observe how, indeed, the presence of the
truncated sum for $\zeta (3)$ in {\eqref{eq:apery3:2}} makes the limit
{\eqref{eq:apery3:lim}} transparent. While, nowadays, it is routine
\cite{schneider-apery} to verify that the sum {\eqref{eq:apery3:2}}
satisfies the recurrence {\eqref{eq:apery3:rec}}, it is much less clear how to
discover sums like {\eqref{eq:apery3:2}} that are suitable for proving a
desired Ap\'ery limit.

Shortly after, and inspired by, Ap\'ery's proof, Beukers
\cite{beukers-irr-leg} gave shorter and more elegant proofs of the
irrationality of $\zeta (2)$ and $\zeta (3)$ by considering double and triple
integrals that result in small linear forms in the zeta value and $1$. For
instance, for $\zeta (3)$, Beukers establishes that
\begin{align}
  & (- 1)^n \int_0^1 \int_0^1 \int_0^1 \frac{x^n (1 - x)^n y^n (1 - y)^n z^n (1
  - z)^n}{(1 - (1 - x y) z)^{n + 1}} \md x \md y \md z \nonumber\\ ={} & A (n) \zeta
  (3) - 6 B (n), \label{eq:linearform:beukers:zeta3}
\end{align}
where $A (n)$ and $B (n)$ are precisely the Ap\'ery numbers
{\eqref{eq:apery3}} and the corresponding secondary solution
{\eqref{eq:apery3:2}}. By bounding the integrand, it is straightforward to
show that the triple integral approaches $0$ as $n \rightarrow \infty$. From
this we directly obtain the Ap\'ery limit {\eqref{eq:apery3:lim}}, namely,
$\lim_{n \rightarrow \infty} B (n) / A (n) = \zeta (3) / 6$.

\begin{example}
  In the same spirit, the Ap\'ery limit {\eqref{eq:delannoy:x:L}} can be
  deduced from
  \begin{equation*}
    \int_0^1 \frac{t^n (1 - t)^n}{(x + t)^{n + 1}} \md t = A_n (x) \ln
     \left(1 + \frac{1}{x} \right) - 2 B_n (x),
  \end{equation*}
  which holds for $x > 0$ with $A_n (x)$ and $B_n (x)$ as in
  Section~\ref{sec:delannoy}. We note that this integral is a variation of the
  integral that Alladi and Robinson \cite{ar-irr} have used to prove
  explicit irrationality measures for numbers of the form $\ln (1 + \lambda)$
  for certain algebraic $\lambda$.
\end{example}

As a powerful variation of this approach, the same kind of linear forms can be
constructed through hypergeometric sums obtained from rational functions. For
instance, Zudilin \cite{zudilin-arithmetic-odd} studies a general
construction, a special case of which is the relation, originally due to
Gutnik,
\begin{equation*}
  - \frac{1}{2} \sum_{t = 1}^{\infty} R_n' (t) = A (n) \zeta (3) - 6 B (n),
   \quad \text{where } R_n (t) = \left(\frac{(t - 1) \cdots (t - n)}{t (t +
   1) \cdots (t + n)} \right)^2,
\end{equation*}
which once more equals {\eqref{eq:linearform:beukers:zeta3}}. We refer to
\cite{zudilin-arithmetic-odd}, \cite{zudilin-appr-polylog} and
\cite[Section~2.3]{avz-apery-limits} for further details and references. A
detailed discussion of the case of $\zeta (2)$ is included in
\cite[Sections~9.5 and 9.6]{bvsz-cf}.

Beukers \cite{beukers-irr} further realized that, in Ap\'ery's cases, the
differential equations associated to the recurrences have a description in
terms of modular forms. Zagier \cite{zagier4} has used such modular
parametrizations to prove Ap\'ery limits in several cases, including for the
Franel numbers, the case $d = 3$ in {\eqref{eq:franel:d}}. The limits occuring
in his cases are rational multiples of
\begin{equation*}
  \zeta (2), \quad L_{- 3} (2) = \sum_{n = 1}^{\infty} \frac{\left(\frac{-
   3}{n} \right)}{n^2}, \quad L_{- 4} (2) = \sum_{n = 1}^{\infty} \frac{\left(\frac{- 4}{n} \right)}{n^2} = \sum_{n = 0}^{\infty} \frac{(- 1)^n}{(2 n +
   1)^2},
\end{equation*}
where $\left(\frac{a}{n} \right)$ is a Legendre symbol and $L_{- 4} (2)$ is
Catalan's constant (whose irrationality remains famously unproven). A general
method for obtaining Ap\'ery limits in cases of modular origin has been
described by Yang \cite{yang-apery-limits}, who proves various Ap\'ery
limits in terms of special values of $L$-functions.

\section{Sums of powers of binomials}\label{sec:franel:d}

Let us consider the family
\begin{equation}
  A^{(d)} (n) = \sum_{k = 0}^n \binom{n}{k}^d \label{eq:franel:d}
\end{equation}
of sums of powers of binomial coefficients. It is easy to see that $A^{(1)}
(n) = 2^n$ and $A^{(2)} (n) = \binom{2 n}{n}$. The numbers $A^{(3)} (n)$ are
known as {\emph{Franel numbers}} \cite[A000172]{oeis}. Almost a century
before the availability of computer-algebra approaches like creative
telescoping, Franel \cite{franel94} obtained explicit recurrences for
$A^{(3)} (n)$ as well as, in a second paper, $A^{(4)} (n)$, and he conjectured
that $A^{(d)} (n)$ satisfies a linear recurrence of order $\lfloor (d + 1) / 2
\rfloor$ with polynomial coefficients. This conjecture was proved by Stoll in
\cite{stoll-rec-bounds}, to which we refer for details and references. It
remains an open problem to show that, in general, no recurrence of lower order
exists.

Van der Poorten \cite[p.~202]{alf} reports that the following Ap\'ery
limits in the cases $d = 3$ and $d = 4$ (in which case the binomial sums
satisfy second-order recurrences like Ap\'ery's sequences) have been
numerically observed by Tom Cusick:
\begin{equation}
  \lim_{n \rightarrow \infty} \frac{B^{(3)} (n)}{A^{(3)} (n)} =
  \frac{\pi^2}{24}, \quad \lim_{n \rightarrow \infty} \frac{B^{(4)}
  (n)}{A^{(4)} (n)} = \frac{\pi^2}{30} . \label{eq:franel:L34}
\end{equation}
In each case, $B^{(d)} (n)$ is the (unique) secondary solution with initial
conditions $B^{(d)} (0) = 0$, $B^{(d)} (1) = 1$. The case $d = 3$ was proved
by Zagier \cite{zagier4} using modular forms. Since the case $d = 4$ is
similarly connected to modular forms \cite{cooper-level10}, we expect that
it can be established using the methods in \cite{yang-apery-limits},
\cite{zagier4} as well. Based on numerical evidence, following the approach
in Section~\ref{sec:search}, we make the following general conjecture
extending {\eqref{eq:franel:L34}}:

\begin{conjecture}
  \label{conj:franel:2}For $d \geq 3$, the minimal-order recurrence
  satisfied by $A^{(d)} (n)$ has a unique solution $B^{(d)} (n)$ with $B^{(d)}
  (0) = 0$ and $B^{(d)} (1) = 1$ that also satisfies
  \begin{equation}
    \lim_{n \rightarrow \infty} \frac{B^{(d)} (n)}{A^{(d)} (n)} = \frac{\zeta
    (2)}{d + 1} . \label{eq:conj:franel:2}
  \end{equation}
\end{conjecture}

Note that for $d \geq 5$, the recurrence is of order $\geq 3$, and
so the solution $B^{(d)} (n)$ is not uniquely determined by the two initial
conditions $B^{(d)} (0) = 0$ and $B^{(d)} (1) = 1$.
Conjecture~\ref{conj:franel:2} asserts that precisely one of these solutions
satisfies {\eqref{eq:conj:franel:2}}.

Subsequent to making this conjecture, we realized that the case $d = 5$ was
already conjectured in \cite[Section~4.1]{avz-apery-limits} as sequence
\#22. We are not aware of previous conjectures for the cases $d \geq 6$.
We have numerically confirmed each of the cases $d \leq 10$ to more than
$100$ digits.

\begin{example}
  \label{eg:franel:5}For $d = 5$, the sum $A^{(5)} (n)$ satisfies a recurrence
  of order $3$, spelled out in \cite{perlstadt-franel}, of the form
  \begin{equation}
    (n + 3)^4 p (n + 1) u (n + 3) + \ldots + 32 (n + 1)^4 p (n + 2) u (n) = 0
    \label{eq:franel:5:rec}
  \end{equation}
  where $p (n) = 55 n^2 + 33 n + 6$. The solution $B^{(5)} (n)$ from
  Conjecture~\ref{conj:franel:2} is characterized by $B^{(5)} (0) = 0$ and
  $B^{(5)} (1) = 1$ and insisting that the recurrence
  {\eqref{eq:franel:5:rec}} also holds for $n = - 1$ (note that this does not
  require a value for $B^{(5)} (- 1)$ because of the term $(n + 1)^4$).
  Similarly, for $d = 6, 7, 8, 9$ the sequence $B^{(d)} (n)$ in
  Conjecture~\ref{conj:franel:2} can be characterized by enforcing the
  recurrence for small negative $n$ and by setting $B^{(d)} (n) = 0$ for $n <
  0$. By contrast, for $d = 10$, there is a two-dimensional space of sequences
  $u (n)$ solving {\eqref{eq:franel:5:rec}} for all integers $n$ with the
  constraint that $u (n) = 0$ for $n \leq 0$. Among these, $B^{(10)} (n)$
  is characterized by $B^{(10)} (1) = 1$ and $B^{(10)} (2) = 381 / 4$.
\end{example}

Now return to the case $d = 5$ and let $C^{(5)} (n)$ be the third solution to
the same recurrence with $C^{(5)} (0) = 0$, $C^{(5)} (1) = 1$, $C^{(5)} (2) =
\frac{48}{7}$. Numerical evidence suggests that we have the Ap\'ery limits
\begin{equation*}
  \lim_{n \rightarrow \infty} \frac{B^{(5)} (n)}{A^{(5)} (n)} = \frac{1}{6}
   \zeta (2), \quad \lim_{n \rightarrow \infty} \frac{C^{(5)} (n)}{A^{(5)}
   (n)} = \frac{3 \pi^4}{1120} = \frac{27}{112} \zeta (4) .
\end{equation*}
Extending this idea to $d = 5, 6, \ldots, 10$, we numerically find Ap\'ery
limits $C^{(d)} (n) / A^{(d)} (n) \rightarrow \lambda \zeta (4)$ with the
following rational values for $\lambda$:
\begin{equation*}
  \frac{27}{112}, \frac{4}{21}, \frac{37}{240}, \frac{7}{55}, \frac{47}{440},
   \frac{1}{11} .
\end{equation*}
These values suggest that $\lambda$ can be expressed as a simple rational
function of $d$:

\begin{conjecture}
  For $d \geq 5$, the minimal-order recurrence satisfied by $A^{(d)} (n)$
  has a unique solution $C^{(d)} (n)$ with $C^{(d)} (0) = 0$ and $C^{(d)} (1)
  = 1$ that also satisfies
  \begin{equation*}
    \lim_{n \rightarrow \infty} \frac{C^{(d)} (n)}{A^{(d)} (n)} = \frac{3 (5
     d + 2)}{(d + 1) (d + 2) (d + 3)} \zeta (4) .
  \end{equation*}
\end{conjecture}

More generally, we expect that, for $d \geq 2 m + 1$, there exist such
Ap\'ery limits for rational multiples of $\zeta (2), \zeta (4), \ldots,
\zeta (2 m)$. It is part of the challenge presented here to explicitly
describe all of these limits. As communicated to us by Zudilin, one could
approach the challenge, uniformly in $d$, by considering the rational
functions
\begin{equation*}
  R_n^{(d)} (t) = \left(\frac{(- 1)^t n!}{t (t + 1) \cdots (t + n)}
   \right)^d
\end{equation*}
in the spirit of \cite{zudilin-arithmetic-odd},
\cite{zudilin-appr-polylog} and \cite[Section~2.3]{avz-apery-limits}, as
indicated in Section~\ref{sec:pf}.

\section{Further challenges and open problems}\label{sec:open}

Although many things are known about Ap\'ery limits, much deserves to be
better understood. The explicit conjectures we offer in the previous section
can be supplemented with similar ones for other families of binomial sums. In
addition, many conjectural Ap\'ery limits that were discovered numerically
are listed in \cite[Section~4]{avz-apery-limits} for sequences that arise
from fourth- and fifth-order differential equations of Calabi--Yau type. As
mentioned in Example~\ref{eg:cy}, it would be of particular interest to
establish all these Ap\'ery limits in a uniform fashion.

It is natural to wonder whether a better understanding of Ap\'ery limits can
lead to new irrationality results. Despite considerable efforts and progress
(we refer the reader to \cite{zudilin-arithmetic-odd} and
\cite{brown-apery} as well as the references therein), it remains a
wide-open challenge to prove the irrationality of, say, $\zeta (5)$ or
Catalan's constant. As a recent promising construction in this direction, we
mention Brown's cellular integrals \cite{brown-apery} which are linear forms
in (multiple) zeta values that are constructed to have certain vanishing
properties that make them amenable to irrationality proofs. In particular,
Brown's general construction includes Ap\'ery's results as (unique) special
cases occuring as initial instances.

In another direction, it would be of interest to systematically study
$q$-analogs and, in particular, to generalize from difference equations to
$q$-difference equations. For instance, Amdeberhan and Zeilberger
\cite{az-qapery} offer an Ap\'ery-style proof of the irrationality of the
$q$-analog of $\ln (2)$ based on a $q$-version of the Delannoy numbers
{\eqref{eq:delannoy}}.

Perron's theorem, which we have mentioned after Poincar\'e's
Theorem~\ref{thm:poincare}, guarantees that, for each characteristic root
$\lambda$ of an appropriate difference equation, there exists a solution $u
(n)$ such that $u (n + 1) / u (n)$ approaches $\lambda$. We note that, for
instance, Ap\'ery's linear form {\eqref{eq:linearform:beukers:zeta3}} is
precisely the unique (up to a constant multiple) solution corresponding to the
$\lambda$ of smallest modulus. General tools to explicitly construct such
Perron solutions from the difference equation would be of great utility.

\begin{acknowledgment}
We are grateful to Alan Sokal for improving the
exposition by kindly sharing lots of careful suggestions and comments. We also
thank Wadim Zudilin for helpful comments, including his suggestion at the end
of Section~\ref{sec:franel:d}, and references.
\end{acknowledgment}

\end{document}